\documentclass{amsart}
\usepackage{amsmath}
\usepackage{amssymb}
\usepackage{amsfonts}
\usepackage{graphicx}

\newtheorem{theorem}{Theorem}
\theoremstyle{plain}

\newtheorem{lemma}{Lemma}

\newtheorem{proposition}{Proposition}

\numberwithin{equation}{section}
\input{tcilatex}

\begin{document}
\title[The Class Averaging Problem]{Representation theoretic patterns in
three-dimensional Cryo-Electron Microscopy II - The class averaging problem}
\author{Ronny Hadani}
\curraddr{Department of Mathematics, University of Texas at Austin, Austin
C1200, USA. }
\author{Amit Singer}
\address{Department of Mathematics and PACM, Princeton University, Fine
Hall, Washington Road, Princeton NJ 08544-1000, USA}
\email{hadani@math.utexas.edu}
\email{amits@math.princeton.edu}
\date{April, 2011}

\begin{abstract}
In this paper we study the formal algebraic structure underlying \ the
intrinsic classification algorithm, recently introduced by Hadani, Shkolnisky, Singer and Zhao, for
classifying noisy projection images of similar viewing directions in
three-dimensional cryo-electron microscopy (cryo-EM). This preliminary
classification is of fundamental importance in determining the
three-dimensional structure of macromolecules from cryo-EM images.
Inspecting this algebraic structure we obtain a conceptual explanation for
the admissibility (correctness) of the algorithm and a proof of its
numerical stability, thus putting it on firm mathematical grounds.
The proof relies on studying the spectral properties of
an integral operator of geometric origin on the two-dimensional sphere,
called the localized parallel transport operator. Along the way, we continue
to develop the representation theoretic setup for three-dimensional cryo-EM
that was initiated in \cite{R}.
\end{abstract}

\maketitle

\section{Introduction}

The goal in cryo-EM is to determine the three-dimensional structure of a
molecule from noisy projection images taken at unknown random orientations
by an electron microscope, i.e., a random Computational Tomography (CT).
Determining three-dimensional structures of large biological molecules
remains vitally important, as witnessed, for example, by the 2003 Chemistry
Nobel Prize, co-awarded to R. MacKinnon for resolving the three-dimensional
structure of the Shaker K+ channel protein \cite{D,M}, and by the 2009
Chemistry Nobel Prize, awarded to V. Ramakrishnan, T. Steitz and A. Yonath
for studies of the structure and function of the ribosome. The standard
procedure for structure determination of large molecules is X-ray
crystallography. The challenge in this method is often more in the
crystallization itself than in the interpretation of the X-ray results,
since many large molecules, including various types of proteins have so far
withstood all attempts to crystallize them.

\begin{figure}[h]
 \centering
 \includegraphics[width=1.8498in,height=2.4197in]{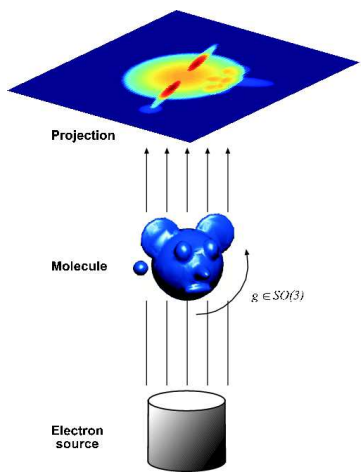}
 % mickymodel.eps: 0x0 pixel, 300dpi, 0.00x0.00 cm, bb=0 0 402 527
 \caption{Schematic drawing of the imaging process: every projection image corresponds to some unknown spatial orientation of the molecule.}
 \label{Mickey_fig}
\end{figure}

Cryo-EM is an alternative approach to X-ray crystallography. In this
approach, samples of identical molecules are rapidly immobilized in a thin
layer of vitreous ice (this is an ice without crystals). The cryo-EM imaging
process produces a large collection of tomographic projections,
corresponding to many copies of the same molecule, each immobilized in a
different (yet unknown) orientation. The intensity of the pixels in a given
projection image is correlated with the line integrals of the electric
potential induced by the molecule along the path of the imaging electrons
(see Figure \ref{Mickey_fig}). The goal is to reconstruct the
three-dimensional structure of the molecule from such a collection of
projection images. The main problem is that the highly intense electron beam
damages the molecule and, therefore, it is problematic to take projection
images of the same molecule at known different directions as in the case of
classical CT\footnote{%
We remark that there are other methods like single-or multi-axis tilt EM
tomogrophy, where several lower dose/higher noise images of a single
molecule are taken from known directions. These methods are used for example
when one has an organic object in vitro or a collection of different objects
in the sample. There is a rich literature for this field starting with the
work of Crowther, DeRosier and Klug in the early 1960s.}. In other words, a
single molecule is imaged only once, rendering an extremely low
signal-to-noise ratio (SNR), mostly due to shot noise induced by the maximal
allowed electron dose.

%\FRAME{ftbphFU}{1.8498in}{2.4197in}{0pt}{\Qcb{%
%Schematic drawing of the imaging process: every projection image corresponds
%to some unknown spatial orientation of the molecule.}}{\Qlb{Mickey_fig}}{%
%mickymodel.eps}{\special{language "Scientific Word";type
%"GRAPHIC";maintain-aspect-ratio TRUE;display "USEDEF";valid_file "F";width
%1.8498in;height 2.4197in;depth 0pt;original-width 5.2434in;original-height
%6.8753in;cropleft "0";croptop "1";cropright "1";cropbottom "0";filename
%'mickymodel.eps';file-properties "XNPEU";}}

\subsection{Mathematical model}

Instead of thinking of a multitude of molecules immobilized in various
orientations and observed by a microscope held in a fixed position, it is
more convenient to think of a single molecule, observed by an electron
microscope from various orientations. Thus, an orientation describes a
configuration of the microscope instead of that of the molecule.

Let $\left( V,\left( \cdot ,\cdot \right) \right) $ be an oriented
three-dimensional Euclidean vector space. The reader can take $V$ to be $%
%TCIMACRO{\U{211d} }%
%BeginExpansion
\mathbb{R}
%EndExpansion
^{3}$ and $\left( \cdot ,\cdot \right) $ to be the standard inner product.
Let $X=\mathrm{Fr}\left( V\right) $ be the oriented frame manifold
associated to $V$; a point $x\in X$ is an orthonormal basis $x=\left(
e_{1},e_{2},e_{3}\right) $ of $V$ compatible with the orientation. The third
vector $e_{3}$ is distinguished, denoted by $\pi \left( x\right) $ and
called \textit{the viewing direction}. More concretely, if we identify $V$
with $%
%TCIMACRO{\U{211d} }%
%BeginExpansion
\mathbb{R}
%EndExpansion
^{3}$, then a point in $X$ can be thought of as a matrix belonging to the
special orthogonal group $SO\left( 3\right) $, whose first,second and third
columns are the vectors $e_{1},e_{2}$ and $e_{3}$ respectively.

\begin{figure}[h]
 \centering
 \includegraphics[width=2.0781in,height=2.2044in]{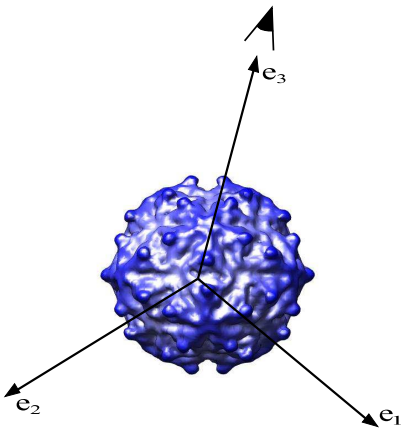}
 % Frame.eps: 0x0 pixel, 300dpi, 0.00x0.00 cm, bb=0 0 446 473
 \caption{A frame $x=\left(
e_{1},e_{2},e_{3}\right) $ modeling the orientation of the electron
microscope, where $\protect\pi \left( x\right) =e_{3}$ is the viewing
direction and the pair $\left( e_{1},e_{2}\right) $ establishes the
coordinates of the camera.}
 \label{Frame_fig}
\end{figure}

Using this terminology, the physics of cryo-EM is modeled as follows:

\begin{itemize}
\item The molecule is modeled by a real valued function $\phi :V\rightarrow 
%TCIMACRO{\U{211d} }%
%BeginExpansion
\mathbb{R}
%EndExpansion
$, describing the electromagnetic potential induced from the charges in the
molecule.

\item A spatial orientation of the microscope is modeled by an orthonormal
frame $x\in X$. The third vector $\pi \left( x\right) $ is the viewing
direction of the microscope and the plane spanned by the first two vectors $%
e_{1}$ and $e_{2}$ is the plane of the camera equipped with the coordinate
system of the camera (see Figure \ref{Frame_fig}).

\item The projection image obtained by the microscope, when observing the
molecule from a spatial orientation $x$ is a real valued function $I:%
%TCIMACRO{\U{211d} }%
%BeginExpansion
\mathbb{R}
%EndExpansion
^{2}\rightarrow 
%TCIMACRO{\U{211d} }%
%BeginExpansion
\mathbb{R}
%EndExpansion
$, given by the X-ray projection along the viewing direction: 
\begin{equation*}
I\left( p,q\right) =\mathrm{Xray}_{\pi \left( x\right) }\phi \left(
p,q\right) =\dint\limits_{t\in 
%TCIMACRO{\U{211d} }%
%BeginExpansion
\mathbb{R}
%EndExpansion
}\phi \left( pe_{1}+qe_{2}+te_{3}\right) dr\text{.}
\end{equation*}
\end{itemize}

for every $\left( p,q\right) \in 
%TCIMACRO{\U{211d} }%
%BeginExpansion
\mathbb{R}
%EndExpansion
^{2}$.

The data collected from the experiment is a set consisting of $N$ projection
images $\mathcal{P=\{}I_{1},..,I_{N}\}$. Assuming that the potential
function $\phi $ is generic\footnote{%
This assumption about the potential $\phi $ can be omitted in the context of
the class averaging algorithm presented in this paper. In particular, the
algorithm can be applied to potentials describing molecules with symmetries
which do not satisfy the "generic" assumption.}, in the sense that, each
image $I_{i}\in \mathcal{P}$ can originate from a unique frame $x_{i}\in X$,
the main problem of cryo-EM is to reconstruct the (unique) unknown frame $%
x_{i}\in X$ associated with each projection image $I_{i}\in \mathcal{P}.$

%\FRAME{ftbpFU}{2.0781in}{2.2044in}{0pt}{\Qcb{A frame $x=\left(
%e_{1},e_{2},e_{3}\right) $ modeling the orientation of the electron
%microscope, where $\protect\pi \left( x\right) =e_{3}$ is the viewing
%direction and the pair $\left( e_{1},e_{2}\right) $ establishes the
%coordinates of the camera.}}{\Qlb{Frame_fig}}{frame.eps}{\special{language
%"Scientific Word";type "GRAPHIC";maintain-aspect-ratio TRUE;display
%"USEDEF";valid_file "F";width 2.0781in;height 2.2044in;depth
%0pt;original-width 6.1886in;original-height 6.5665in;cropleft "0";croptop
%"1";cropright "1";cropbottom "0";filename
%'../../../cryo_project/symmetries/Frame.eps';file-properties "XNPEU";}}

\subsection{Class averaging}

\begin{figure}[h]
 \centering
 \includegraphics[width=4.382in,height=1.2886in]{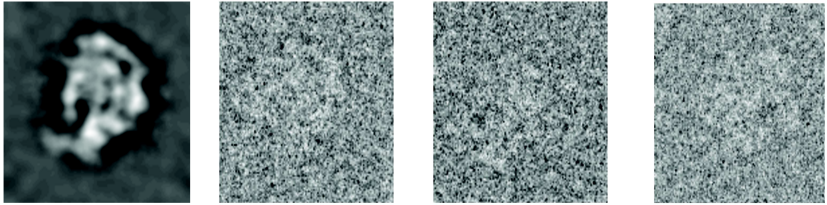}
 % ClassAverage.eps: 0x0 pixel, 300dpi, 0.00x0.00 cm, bb=0 0 493 144
 \caption{The left most image is a clean
simulated projection image of the E.coli 50S ribosomal subunit. The other
three images are real electron microscope images of the same subunit.}
 \label{Class_fig}
\end{figure}

As projection images in cryo-EM have extremely low SNR\footnote{%
SNR stands for Signal to Noise Ratio, which is the ratio between the squared 
$L^{2}$ norm of the signal and the squared $L^{2}$ norm of the noise.} (see
Figure \ref{Class_fig}), a crucial initial step in all reconstruction
methods is \textquotedblleft class averaging\textquotedblright\ \cite{F}.
Class averaging is the grouping of a large data set of noisy raw projection
images into clusters, such that images within a single cluster have similar
viewing directions. Averaging rotationally aligned noisy images within each
cluster results in \textquotedblleft class averages\textquotedblright ;
these are images that enjoy a higher SNR and are used in later cryo-EM
procedures such as the angular reconstitution procedure, \cite{V}, that
requires better quality images. Finding consistent class averages is
challenging due to the high level of noise in the raw images.

%\FRAME{ftbpFU}{4.382in}{1.2886in}{0pt}{\Qcb{The left most image is a clean
%simulated projection image of the E.coli 50S ribosomal subunit. The other
%three images are real electron microscope images of the same subunit.}}{\Qlb{%
%Class_fig}}{classaverage.eps}{\special{language "Scientific Word";type
%"GRAPHIC";maintain-aspect-ratio TRUE;display "USEDEF";valid_file "F";width
%4.382in;height 1.2886in;depth 0pt;original-width 6.8165in;original-height
%1.9856in;cropleft "0";croptop "1";cropright "1";cropbottom "0";filename
%'ClassAverage.eps';file-properties "XNPEU";}}

\begin{center}
\medskip
\end{center}

The starting point for the classification is the idea that visual similarity
between projection images suggests vicinity between viewing directions of
the corresponding (unknown) frames. The similarity between images $I_{i}$
and $I_{j}$ is measured by their invariant distance (introduced in \cite{P})
which is the Euclidean distance between the images when they are optimally
aligned with respect to in-plane rotations, namely%
\begin{equation}
d\left( I_{i},I_{j}\right) =\min_{g\in SO\left( 2\right) }\left\Vert R\left(
g\right) I_{i}-I_{j}\right\Vert ,  \label{invdist_eq}
\end{equation}%
where 
\begin{equation*}
R\left( g\right) I\left( p,q\right) =I\left( g^{-1}\left( p,q\right) \right)
,
\end{equation*}%
for any function $I:%
%TCIMACRO{\U{211d} }%
%BeginExpansion
\mathbb{R}
%EndExpansion
^{2}\rightarrow 
%TCIMACRO{\U{211d} }%
%BeginExpansion
\mathbb{R}
%EndExpansion
$.

One can choose some threshold value $\epsilon $, such that $d\left(
I_{i},I_{j}\right) \leq \epsilon $ is indicative that perhaps the
corresponding frames $x_{i}$ and $x_{j}$ have nearby viewing directions. The
threshold $\epsilon $ defines an undirected graph $G=\left( \mathrm{Vertices}%
,\mathrm{Edges}\right) $ with vertices labeled by numbers $1,..,N$ and an
edge connecting vertex $i$ with vertex $j$ if and only if the invariant
distance between the projection images $I_{i}$ and $I_{j}$ is smaller then $%
\epsilon $, namely%
\begin{equation*}
\left\{ i,j\right\} \in \mathrm{Edges}\Longleftrightarrow d\left(
I_{i},I_{j}\right) \leq \epsilon \text{.}
\end{equation*}

In an ideal noiseless world, the graph $G$ acquires the geometry of the unit
sphere $S\left( V\right) $, namely, two images are connected by an edge if
and only if their corresponding viewing directions are close on the sphere,
in the sense that they belong to some small spherical cap of opening angle $%
a=a\left( \epsilon \right) $.

However, the real world is far from ideal as it is governed by noise; hence,
it often happens that two images of completely different viewing directions
have small invariant distance. This can happen when the realizations of the
noise in the two images match well for some random in-plane rotation,
leading to spurious neighbor identification. Therefore, the na\"{\i}ve
approach of averaging the rotationally aligned nearest neighbor images can
sometimes yield a poor estimate of the true signal in the reference image.

\medskip

\textbf{To summarize: }From this point of view, the main problem is to
distinguish the good edges from the bad ones in the graph $G$, or, in other
words, to distinguish the true neighbors from the false ones (called \textit{%
outliers).} The existence of outliers is the reason why the classification
problem is non-trivial. We emphasize that without excluding the outliers,
averaging rotationally aligned images of small invariant distance (\ref%
{invdist_eq}) yields poor estimate of the true signal, rendering the problem
of three-dimensional reconstruction from cryo-EM images non-feasible. In
this respect, the class averaging problem is of fundamental importance.

\subsection{Main results}

In \cite{S2}, we introduced a novel algorithm, referred to in this paper as
the \textit{intrinsic classification algorithm}, for classifying noisy
projection images of similar viewing directions. The main appealing property
of this new algorithm is its extreme robustness to noise and to presence of
outliers; in addition, it also enjoys efficient time and space complexity.
These properties are explained thoroughly in \cite{S2}, which includes also
a large number of numerical experiments.

In this paper we study the formal algebraic structure that underlies the
intrinsic classification algorithm. Inspecting this algebraic structure we
obtain a conceptual explanation for the admissibility (correctness) of the
algorithm and a proof of its numerical stability, thus putting\ it on firm
mathematical grounds. The proof relies on the study of a certain integral
operator $T_{h}$ on $X$, of geometric origin, called the \textit{localized
parallel transport operator}. Specifically:

\begin{itemize}
\item Admissibility amounts to the fact that the maximal eigenspace of $%
T_{h} $ is a three-dimensional complex Hermitian vector space and that there
is a \textbf{canonical} identification of Hermitian vector spaces between
this eigenspace and the complexified vector space $W=%
%TCIMACRO{\U{2102} }%
%BeginExpansion
\mathbb{C}
%EndExpansion
V$.

\item Numerical stability amounts to the existence of a spectral gap which
separates the maximal eigenvalue of $T_{h}$ from the rest of the spectrum,
which enables one to obtain a stable numerical approximation of the
corresponding maximal eigenspace and of other related geometric structures.
\end{itemize}

The main technical result of this paper is a complete description of the
spectral properties of the localized parallel transport operator. Along the
way, we continue to develop the mathematical set-up for cryo-EM that was
initiated in \cite{R}, thus further elucidating the central role played by
representation theoretic principles in this scientific discipline.

The remainder of the introduction is devoted to a detailed description of
the intrinsic classification algorithm and to an explanation of the main
ideas and results of this paper.

\subsection{Transport data\label{transport_sub}}

A preliminary step is to extract certain geometric data from the set of
projection images, called (local) \textit{empirical} \textit{transport data}.

When computing the invariant distance between images $I_{i}$ and $I_{j}$ we
also record the rotation matrix in $SO\left( 2\right) $ that realizes the
minimum in (\ref{invdist_eq}) and denote this special rotation by $%
\widetilde{T}\left( i,j\right) $, that is 
\begin{equation}
\widetilde{T}\left( i,j\right) =\underset{g\in SO(2)}{\mathrm{argmin}}%
\left\Vert R\left( g\right) I_{i}-I_{j}\right\Vert .  \label{emp_eq}
\end{equation}%
noting that, 
\begin{equation}
\widetilde{T}\left( j,i\right) =\widetilde{T}\left( i,j\right) ^{-1}.
\label{symm_eq}
\end{equation}

The main observation is that in an ideal noiseless world the rotation $%
\widetilde{T}\left( i,j\right) $ can be interpreted as a geometric relation
between the corresponding frames $x_{i}$ and $x_{j}$, provided the invariant
distance between the corresponding images is small. This relation is
expressed in terms of parallel transport on the sphere, as follows: define
the rotation

\begin{equation*}
T\left( x_{i},x_{j}\right) =%
\begin{pmatrix}
\cos \left( \theta _{ij}\right) & -\sin \left( \theta _{ij}\right) \\ 
\sin \left( \theta _{ij}\right) & \text{ \ }\cos \left( \theta _{ij}\right)%
\end{pmatrix}%
,
\end{equation*}%
as the unique solution of the equation 
\begin{equation}
x_{i}\vartriangleleft T\left( x_{i},x_{j}\right) =t_{\pi \left( x_{i}\right)
,\pi \left( x_{j}\right) }x_{j},  \label{geom_eq}
\end{equation}%
where $t_{\pi \left( x_{i}\right) ,\pi \left( x_{j}\right) }$ is the
parallel transport along the unique geodesic on the sphere connecting the
points $\pi \left( x_{j}\right) $ with $\pi \left( x_{i}\right) $ or, in
other words, it is the rotation in $SO\left( V\right) $ that takes the
vector $\pi \left( x_{j}\right) $ to $\pi \left( x_{i}\right) $ along the
shortest path on the sphere and the action $\vartriangleleft $ is defined by 
\begin{equation*}
x\vartriangleleft 
\begin{pmatrix}
\cos \left( \theta \right) & -\sin \left( \theta \right) \\ 
\sin \left( \theta \right) & \text{ \ }\cos \left( \theta \right)%
\end{pmatrix}%
=(\cos \left( \theta \right) e_{1}+\sin \left( \theta \right) e_{2},-\sin
\left( \theta \right) e_{1}+\cos \left( \theta \right) e_{2},e_{3}),
\end{equation*}%
for every $x=\left( e_{1},e_{2},e_{3}\right) $. The precise statement is
that the rotation $\widetilde{T}\left( i,j\right) $ approximates the
rotation $T\left( x_{i},x_{j}\right) $ when $\left\{ i,j\right\} \in \mathrm{%
Edges}$. This geometric interpretation of the rotation $\widetilde{T}\left(
i,j\right) $ is suggested from\ a combination of mathematical and empirical
considerations that we proceed to explain.

\begin{itemize}
\item On the mathematical side: the rotation $T\left( x_{i},x_{j}\right) $
is the unique rotation of the frame $x_{i}$ around its viewing direction $%
\pi \left( x_{i}\right) $, minimizing the distance to the frame $x_{j}$.
This is a standard fact from differential geometry (a direct proof of this
statement appears in \cite{S2}).

\item On the empirical side: if the function $\phi $ is "nice", then the
optimal alignment $\widetilde{T}\left( i,j\right) $ of the projection images
is correlated with the optimal alignment $T\left( x_{i},x_{j}\right) $ of
the corresponding frames. This correlation of course improves as the
distance between $\pi \left( x_{i}\right) $ and $\pi \left( x_{j}\right) $
becomes smaller. A quantitative study of the relation between $\widetilde{T}%
\left( i,j\right) $ and $T\left( x_{i},x_{j}\right) $ involves
considerations from image processing thus it is beyond the scope of this
paper.
\end{itemize}

To conclude, the "empirical" rotation $\widetilde{T}\left( i,j\right) $
approximates the "geometric" rotation $T\left( x_{i},x_{j}\right) $ only
when the viewing directions $\pi \left( x_{i}\right) $ and $\pi \left(
x_{j}\right) $ are close, in the sense that they belong to some small
spherical cap of opening angle $a$. The later "geometric" condition is
correlated with the "empirical" condition that the corresponding images $%
I_{i}$ and $I_{j}$ have small invariant distance. When $\pi \left(
x_{i}\right) $ and $\pi \left( x_{j}\right) $ are far from each other, the
rotation $\widetilde{T}\left( i,j\right) $ is not related any longer to
parallel transportation on the sphere. For this reason, we consider only
rotations $\widetilde{T}\left( i,j\right) $ for\ which $\left\{ i,j\right\}
\in \mathrm{Edges}$ and call this collection the (local) \textit{empirical} 
\textit{transport data}.

\subsection{The intrinsic classification algorithm}

The intrinsic classification algorithm accepts as an input the empirical
transport data $\{\widetilde{T}\left( i,j\right) :\left\{ i,j\right\} \in 
\mathrm{Edges}\}$ and produces as an output the Euclidean inner products $%
\{(\pi \left( x_{i}\right) ,\pi \left( x_{j}\right) ):i,j=1,..N\}$. Using
these inner products, one can identify the true neighbors in the graph $G$,
as the pairs $\left\{ i,j\right\} \in \mathrm{Edges}$ for which the inner
product $\left( \pi \left( x_{i}\right) ,\pi \left( x_{j}\right) \right) $
is close to $1$. The formal justification of the algorithm requires the
empirical assumption that the frames $x_{i}$, $i=1,..,N$ are uniformly
distributed in the\ frame manifold $X$, according to the unique normalized
Haar measure on $X$. This assumption corresponds to the situation where the
orientations of the molecules in the ice are distributed independently and
uniformly at random.

The main idea of the algorithm is to construct an intrinsic model, denoted
by $\mathbb{W}_{N}$, of the Hermitian vector space $W=%
%TCIMACRO{\U{2102} }%
%BeginExpansion
\mathbb{C}
%EndExpansion
V$ which is expressed solely in terms of the empirical transport data.

The algorithm proceeds as follows:

\medskip

\textbf{Step1 (Ambient Hilbert space):} consider the standard $N$%
-dimensional Hilbert space 
\begin{equation*}
\mathcal{H}_{N}=%
%TCIMACRO{\U{2102} }%
%BeginExpansion
\mathbb{C}
%EndExpansion
^{N}\text{.}
\end{equation*}

\textbf{Step 2 (Self adjoint operator): }identify $%
%TCIMACRO{\U{211d} }%
%BeginExpansion
\mathbb{R}
%EndExpansion
^{2}$ with $%
%TCIMACRO{\U{2102} }%
%BeginExpansion
\mathbb{C}
%EndExpansion
$ and consider each rotation $\widetilde{T}(i,j)$ as a complex number of
unit norm. Define the $N\times N$ complex matrix 
\begin{equation*}
\widetilde{T}_{N}:\mathcal{H}_{N}\rightarrow \mathcal{H}_{N},
\end{equation*}%
by putting the rotation $\widetilde{T}(i,j)$ in the $\left( i,j\right) $
entry. Notice that the matrix $\widetilde{T}_{N}$ is self-adjoint by (\ref%
{symm_eq}).

\textbf{Step 3 (Intrinsic model): }the matrix $\widetilde{T}_{N}$ induces a
spectral decomposition 
\begin{equation*}
\mathcal{H}_{N}\mathcal{=}\bigoplus\limits_{\lambda }\mathcal{H}_{N}\left(
\lambda \right) .
\end{equation*}

\begin{theorem}
\label{spec_thm}There exists a threshold $\lambda _{0}$ such that 
\begin{equation*}
\dim \bigoplus\limits_{\lambda >\lambda _{0}}\mathcal{H}_{N}\left( \lambda
\right) =3\text{.}
\end{equation*}
\end{theorem}

Define the Hermitian vector space 
\begin{equation*}
\mathbb{W}_{N}=\bigoplus\limits_{\lambda >\lambda _{0}}\mathcal{H}_{N}\left(
\lambda \right) \text{.}
\end{equation*}

\textbf{Step 4 (Computation of the Euclidean inner products): }the Euclidean
inner products $\{(\pi \left( x_{i}\right) ,\pi \left( x_{j}\right)
):i,j=1,..N\}$ are computed from the vector space $\mathbb{W}_{N}$, as
follows: for every $i=1,..,N$, denote by $\varphi _{i}\in \mathbb{W}_{N}$
the vector 
\begin{equation*}
\varphi _{i}=\sqrt{2/3}\cdot \mathrm{pr}_{i}^{\ast }\left( 1\right) \text{,}
\end{equation*}%
where $\mathrm{pr}_{i}:\mathbb{W}_{N}\rightarrow 
%TCIMACRO{\U{2102} }%
%BeginExpansion
\mathbb{C}
%EndExpansion
$ is the projection on the $i$th component and $\mathrm{pr}_{i}^{\ast }:%
%TCIMACRO{\U{2102} }%
%BeginExpansion
\mathbb{C}
%EndExpansion
\rightarrow \mathbb{W}_{N}$ is the adjoint map. In addition, for every frame 
$x\in X$, $x=\left( e_{1},e_{2},e_{3}\right) $, denote by $\delta _{x}\in W$
the (complex) vector $e_{1}-ie_{2}$.

The upshot is that the intrinsic vector space $\mathbb{W}_{N}$ consisting of
the collection of vectors $\varphi _{i}\in \mathbb{W}_{N}$, $i=1,..,N$ is
(approximately\footnote{%
This approximation improves as $N$ grows.}) isomorphic to the extrinsic
vector space $W$ consisting of the collection of vectors $\delta _{x_{i}}\in
W$, $i=1,..,N$, where $x_{i}$ is the frame corresponding to the image $I_{i}$%
, for every $i=1,..,N$. This statement is the content of the following
theorem:

\begin{theorem}
\label{isom_thm}There exists a unique (approximated) isomorphism $\tau _{N}:W%
\overset{\simeq }{\rightarrow }\mathbb{W}_{N}$ of Hermitian vector spaces
such that 
\begin{equation*}
\tau _{N}\left( \delta _{x_{i}}\right) =\varphi _{i},
\end{equation*}%
for every $i=1,..,N$.
\end{theorem}

The above theorem enables us to express, in intrinsic terms, the Euclidean
inner products between the viewing directions, as follows: starting with the
following identity from linear algebra (that will be proved in the sequel): 
\begin{equation}
\left( \pi \left( x\right) ,\pi \left( y\right) \right) =\left\vert
\left\langle \delta _{x},\delta _{y}\right\rangle \right\vert -1\text{,}
\label{LinAlg_eq}
\end{equation}%
for every pair of frames $x,y\in X$, where $\left( \cdot ,\cdot \right) $ is
the Euclidean product on $V$ and $\left\langle \cdot ,\cdot \right\rangle $
is the Hermitian product on $W=%
%TCIMACRO{\U{2102} }%
%BeginExpansion
\mathbb{C}
%EndExpansion
V$, induced from $\left( \cdot ,\cdot \right) $, given by%
\begin{equation*}
\left\langle u+iv,u^{\prime }+iv^{\prime }\right\rangle =\left( u,v\right)
+\left( v,v^{\prime }\right) -i\left( u,v^{\prime }\right) +i\left(
v,u^{\prime }\right) \text{,}
\end{equation*}%
we obtain the following relation:%
\begin{equation}
\left( \pi \left( x_{i}\right) ,\pi \left( x_{j}\right) \right) =\left\vert
\left\langle \varphi _{i},\varphi _{j}\right\rangle \right\vert -1,
\label{intrinsic_eq}
\end{equation}%
for every $i,j=1,..,N$. In the derivation of Relation (\ref{intrinsic_eq})
from Relation (\ref{LinAlg_eq}) we use Theorem \ref{isom_thm}. Notice that
Relation (\ref{intrinsic_eq}) implies that although we do not know the frame
associated with every projection image, we still are able to compute the
inner product between every pair of such frames from the intrinsic vector
space $\mathbb{W}_{N}$ which, in turns, can be computed from the images.

\subsection{Structure of the paper}

The paper consists of three sections besides the introduction.

\begin{itemize}
\item In Section \ref{prelim_sec}, we begin by introducing the basic
analytic setup which is relevant for the class averaging problem in cryo-EM.
Then, we proceed to formulate the main results of this paper, which are: a
complete description of the spectral properties of the localized parallel
transport operator (Theorem \ref{spectral_thm}), the spectral gap property
(Theorem \ref{gap_thm}) and the admissibility of the intrinsic
classification algorithm (Theorems \ref{char_thm} and \ref{dist_thm}).

\item In Section \ref{spectral_sec}, we prove Theorem \ref{spectral_thm}: in
particular, we develop all the representation theoretic machinery that is
needed for the proof.

\item Finally, in Appendix \ref{proofs_sec}, we give the proofs of all
technical statements which appear in the previous sections.
\end{itemize}

\bigskip

\textbf{Acknowledgement: }\textit{The first author would like to thank
Joseph Bernstein for many helpful discussions concerning the mathematical
aspects of this work. He also thanks Richard Askey for his valuable advice
about Legendre polynomials. The second author is partially supported by
Award Number R01GM090200 from the National Institute of General Medical
Sciences. The content is solely the responsibility of the authors and does
not necessarily represent the official views of the National Institute of
General Medical Sciences or the National Institutes of Health. This work is
part of a project conducted jointly with Shamgar Gurevich, Yoel Shkolnisky
and Fred Sigworth.}

\section{Preliminaries and main results\label{prelim_sec}}

\subsection{Setup}

Let $\left( V,\left( \cdot ,\cdot \right) \right) $ be a three-dimensional,
oriented, Euclidean vector space over $%
%TCIMACRO{\U{211d} }%
%BeginExpansion
\mathbb{R}
%EndExpansion
$. The reader can take $V$ $=$ $%
%TCIMACRO{\U{211d} }%
%BeginExpansion
\mathbb{R}
%EndExpansion
^{3}$ equipped with the standard orientation and $\left( \cdot ,\cdot
\right) $ to be the standard inner product. Let $W=%
%TCIMACRO{\U{2102} }%
%BeginExpansion
\mathbb{C}
%EndExpansion
V$ denote the complexification of $V$. We equip $W$ with the Hermitian
product $\left\langle \cdot ,\cdot \right\rangle :W\times W\rightarrow 
%TCIMACRO{\U{2102} }%
%BeginExpansion
\mathbb{C}
%EndExpansion
$, induced from $\left( \cdot ,\cdot \right) $, given by 
\begin{equation*}
\left\langle u+iv,u^{\prime }+iv^{\prime }\right\rangle =\left( u,v\right)
+\left( v,v^{\prime }\right) -i\left( u,v^{\prime }\right) +i\left(
v,u^{\prime }\right) \text{.}
\end{equation*}

Let $SO\left( V\right) $ denote the group of orthogonal transformations with
respect to the inner product $\left( \cdot ,\cdot \right) $ which preserve
the orientation. Let $S\left( V\right) $ denote the unit sphere in $V$, that
is, $S\left( V\right) =\left\{ v\in V:\left( v,v\right) =1\right\} $. Let $%
X=Fr\left( V\right) $ denote the manifold of oriented orthonormal frames in $%
V$, that is, a point $x\in X$ is an orthonormal basis $x=\left(
e_{1},e_{2},e_{3}\right) $ of $V$ compatible with the orientation.

We consider two commuting group actions on the frame manifold: a left action
of the group $SO\left( V\right) $, given by 
\begin{equation*}
g\vartriangleright \left( e_{1},e_{2},e_{3}\right) =\left(
ge_{1},ge_{2},ge_{3}\right) ,
\end{equation*}%
and a right action of the special orthogonal group $SO(3)$, given by 
\begin{eqnarray*}
\left( e_{1},e_{2},e_{3}\right) \vartriangleleft g
&=&(a_{11}e_{1}+a_{21}e_{2}+a_{31}e_{3}, \\
&&a_{12}e_{1}+a_{22}e_{2}+a_{32}e_{3}, \\
&&a_{13}e_{1}+a_{23}e_{2}+a_{33}e_{3}),
\end{eqnarray*}%
for 
\begin{equation*}
g=%
\begin{pmatrix}
a_{11} & a_{12} & a_{13} \\ 
a_{21} & a_{22} & a_{23} \\ 
a_{31} & a_{32} & a_{33}%
\end{pmatrix}%
.
\end{equation*}

We distinguish the copy of $SO(2)$ inside $SO(3)$ consisting of matrices of
the form 
\begin{equation*}
g=%
\begin{pmatrix}
a_{11} & a_{12} & 0 \\ 
a_{21} & a_{22} & 0 \\ 
0 & 0 & 1%
\end{pmatrix}%
,
\end{equation*}%
and consider $X$ as a principal $SO\left( 2\right) $ bundle over $S\left(
V\right) $ where the fibration map $\pi :X\rightarrow S\left( V\right) $ is
given by $\pi \left( e_{1},e_{2},e_{3}\right) =e_{3}$. We call the vector $%
e_{3}$ the \textit{viewing direction.}

\subsection{The Transport data}

Given a point $v\in S\left( V\right) $, we denote by $X_{v}$ the fiber of
the frame manifold laying over $v$, that is, $X_{v}=\left\{ x\in X:\pi
\left( x\right) =v\right\} $. For every pair of frames $x,y\in X$ such that $%
\pi \left( x\right) \neq \pm \pi \left( y\right) $, we define a matrix $%
T\left( x,y\right) \in SO\left( 2\right) $, characterized by the property 
\begin{equation*}
x\vartriangleleft T\left( x,y\right) =t_{\pi \left( x\right) ,\pi \left(
y\right) }\left( y\right) ,
\end{equation*}%
where $t_{\pi \left( x\right) ,\pi \left( y\right) }:X_{\pi \left( y\right)
}\rightarrow X_{\pi \left( x\right) }$ is the morphism between the
corresponding fibers, given by the parallel transport mapping along the
unique geodesic in the sphere $S\left( V\right) $ connecting the points $\pi
\left( y\right) $ with $\pi \left( x\right) $. We identify $%
%TCIMACRO{\U{211d} }%
%BeginExpansion
\mathbb{R}
%EndExpansion
^{2}$ with $%
%TCIMACRO{\U{2102} }%
%BeginExpansion
\mathbb{C}
%EndExpansion
$ and consider\ $T\left( x,y\right) $ as a complex number of unit norm. The
collection of matrices $\left\{ T\left( x,y\right) \right\} $ satisfy the
following properties:

\begin{itemize}
\item \textbf{Symmetry: }For every $x,y\in X$, we have $T\left( x,y\right)
=T\left( x,y\right) ^{-1}$, where the left hand side of the equality
coincides with the complex conjugate $\overline{T\left( x,y\right) }$. This
property follows from the fact that the parallel transport mapping
satisfies: 
\begin{equation*}
t_{\pi \left( y\right) ,\pi \left( x\right) }=t_{\pi \left( x\right) ,\pi
\left( y\right) }^{-1}.
\end{equation*}

\item \textbf{Invariance: }For every $x,y\in X$ and element $g\in SO\left(
V\right) $, we have that $T\left( g\vartriangleright x,g\vartriangleright
y\right) =T\left( x,y\right) $. This property follows from the fact that the
parallel transport mapping satisfies: 
\begin{equation*}
t_{\pi \left( g\vartriangleright x\right) ,\pi \left( g\vartriangleright
y\right) }=g\circ t_{\pi \left( x\right) ,\pi \left( y\right) }\circ g^{-1},
\end{equation*}%
for every $g\in SO\left( V\right) $.

\item \textbf{Equivariance: }For every $x,y\in X$ and elements $%
g_{1},g_{2}\in SO\left( 2\right) $, we have that $T\left( x\vartriangleleft
g_{1},y\vartriangleleft g_{2}\right) =g_{1}^{-1}T\left( x,y\right) g_{2}.$
This property follows from the fact that the parallel transport mapping
satisfies: 
\begin{equation*}
t_{\pi \left( x\vartriangleleft g_{1}\right) ,\pi \left( y\vartriangleleft
g_{2}\right) }=t_{\pi \left( x\right) ,\pi \left( y\right) },
\end{equation*}
for every $g_{1},g_{2}\in SO\left( 2\right) $.
\end{itemize}

The collection $\left\{ T(x,y)\right\} $ is referred to as the \textit{%
transport data}.

\subsection{The parallel transport operator \label{parallel_sub}}

Let $\mathcal{H=}%
%TCIMACRO{\U{2102} }%
%BeginExpansion
\mathbb{C}
%EndExpansion
\left( X\right) $ denote the Hilbertian space of smooth complex valued
functions on $X$ (here, the word Hilbertian means that $\mathcal{H}$ is not
complete)\footnote{%
In general, in this paper, we will not distinguish between an Hilbertian
vector space and its completion and the correct choice between the two will
be clear from the context.}, where the Hermitian product is the standard
one, given by 
\begin{equation*}
\left\langle f_{1},f_{2}\right\rangle _{\mathcal{H}}=\int\limits_{x\in
X}f_{1}\left( x\right) \overline{f_{2}\left( x\right) }dx,
\end{equation*}%
for every $f_{1},f_{2}\in \mathcal{H}$, where $dx$ denotes the normalized
Haar measure on $X$. In addition, $\mathcal{H}$ supports a unitary
representation of the group $SO\left( V\right) \times SO\left( 2\right) $,
where the action of an element $g=\left( g_{1},g_{2}\right) $ sends a
function $s\in \mathcal{H}$ to a function $g\cdot s$, given by 
\begin{equation*}
\left( g\cdot s\right) \left( x\right) =s\left( g_{1}^{-1}\vartriangleright
x\vartriangleleft g_{2}\right) ,
\end{equation*}%
for every $x\in X$.

Using the transport data, we define an integral operator $T:\mathcal{H}%
\rightarrow \mathcal{H}$ as 
\begin{equation*}
T\left( s\right) \left( x\right) =\int\limits_{y\in X}T\left( x,y\right)
s\left( y\right) dy,
\end{equation*}%
for every $s\in \mathcal{H}$. The properties of the transport data imply the
following properties of the operator $T$:

\begin{itemize}
\item The symmetry property implies that $T$ is self adjoint.

\item The invariance property implies that $T$ commutes with the $SO\left(
V\right) $ action, namely $T\left( g\cdot s\right) =g\cdot T\left( s\right) $
for every $s\in \mathcal{H}$ and $g\in SO\left( V\right) $.

\item The implication of the equivariance property will be discussed later
when we study the kernel of $T$.
\end{itemize}

The operator $T$ is referred to as the \textit{parallel transport operator}.

\subsubsection{Localized parallel transport operator}

The operator which arise naturally in our context is a localized version of
the transport operator. Let us fix an angle $a\in \left[ 0,\pi \right] $,
designating an opening angle of a spherical cap on the sphere and consider
the parameter $h=1-\cos (a)$, taking values in the interval $\left[ 0,2%
\right] $.

Given a choice of this parameter, we define an integral operator $T_{h}:%
\mathcal{H\rightarrow H}$ as 
\begin{equation}
T_{h}\left( s\right) \left( x\right) =\int\limits_{y\in B\left( x,a\right)
}T\left( x,y\right) s\left( y\right) dy.  \label{loc_eq}
\end{equation}%
where $B\left( x,a\right) =\left\{ y\in X:(\pi \left( x\right) ,\pi \left(
y\right) )>\cos \left( a\right) \right\} $. Similar considerations as before
show that $T_{h}$ is self-adjoint and, in addition, commutes with the $%
SO\left( V\right) $ action. Finally, note that the operator $T_{h}$ should
be considered as a localization of the operator of parallel transport
discussed in the previous paragraph, in the sense, that now only frames with
close viewing directions interact through the integral (\ref{loc_eq}). For
this reason, the operator $T_{h}$ is referred to as the \textit{localized
parallel transport operator}.

\subsection{Spectral properties of the localized parallel transport operator}

We focus our attention on the spectral properties of the operator $T_{h}$,
in the regime $h\ll 1$, since this is the relevant regime for the class
averaging application.

\begin{theorem}
\label{spectral_thm}The operator $T_{h}$ has a discrete spectrum $\lambda
_{n}\left( h\right) $, $n\in 
%TCIMACRO{\U{2115} }%
%BeginExpansion
\mathbb{N}
%EndExpansion
$, such that $\dim \mathcal{H}\left( \lambda _{n}\left( h\right) \right)
=2n+1$, for every $h\in (0,2]$, Moreover, in the regime $h\ll 1$, the
eigenvalue $\lambda _{n}\left( h\right) $ has the asymptotic expansion 
\begin{equation*}
\lambda _{n}\left( h\right) =\frac{1}{2}h-\frac{1+\left( n+2\right) \left(
n-1\right) }{8}h^{2}+O\left( h^{3}\right) .
\end{equation*}
\end{theorem}

For a proof, see Section \ref{spectral_sec}.

In fact, each eigenvalue $\lambda _{n}\left( h\right) $, as a function of $h$%
, is a polynomial of degree $n+1$. In Section \ref{spectral_sec}, we give a
complete description of these polynomials by means of a generating function.
To get some feeling for the formulas that arise, we list below the first
four eigenvalues%
\begin{eqnarray*}
\lambda _{1}\left( h\right) &=&\frac{1}{2}h-\frac{1}{8}h^{2}, \\
\lambda _{2}\left( h\right) &=&\frac{1}{2}h-\frac{5}{8}h^{2}+\frac{1}{6}%
h^{3}, \\
\lambda _{3}\left( h\right) &=&\frac{1}{2}h-\frac{11}{8}h^{2}+\frac{25}{24}%
h^{3}-\frac{15}{64}h^{4}, \\
\lambda _{4}\left( h\right) &=&\frac{1}{2}h-\frac{19}{8}h^{2}+\frac{27}{8}%
h^{3}-\frac{119}{64}h^{4}+\frac{7}{20}h^{5}\text{.}
\end{eqnarray*}

\begin{figure}[t]
 \centering
 \includegraphics[width=3.813in,height=2.6247in]{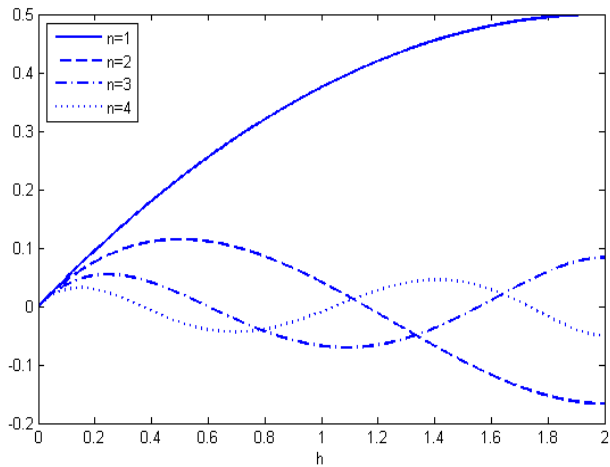}
 % eigenvalues1234.eps: 0x0 pixel, 300dpi, 0.00x0.00 cm, bb=0 0 513 323
 \caption{The first four eigen values of
the operator $T_{h}$, presented as functions of $h\in \lbrack 0,2]$.}
 \label{Eigen1234_fig}
\end{figure}

The graphs of $\lambda _{i}\left( h\right) $, $i=1,2,3,4$ are given in
Figure \ref{Eigen1234_fig}.

%\FRAME{ftbpFU}{3.813in}{2.6247in}{0pt}{\Qcb{The first four eigen values of
%the operator $T_{h}$, presented as functions of $h\in \lbrack 0,2]$. }}{\Qlb{%
%Eigen1234_fig}}{eigenvalues1234.eps}{\special{language "Scientific
%Word";type "GRAPHIC";display "USEDEF";valid_file "F";width 3.813in;height
%2.6247in;depth 0pt;original-width 6.6565in;original-height 4.1563in;cropleft
%"0";croptop "1";cropright "1";cropbottom "0";filename
%'eigenvalues1234.tif';file-properties "XNPEU";}}

\subsubsection{Spectral gap}

Noting that $\lambda _{2}\left( h\right) $ attains its maximum at $h=1/2$,
we have

\begin{theorem}
\label{gap_thm}For every value of $h\in \left[ 0,2\right] $, the maximal
eigenvalue of $T_{h}$ is $\lambda _{1}\left( h\right) $. Moreover, for every
value of $h\in \left[ 0,1/2\right] $, there is a spectral gap $G(h)$ of the
form 
\begin{equation*}
G(h)=\lambda _{1}\left( h\right) -\lambda _{2}\left( h\right) =\frac{1}{2}%
h^{2}-\frac{1}{6}h^{3}.
\end{equation*}
\end{theorem}

For a proof, see Appendix \ref{proofs_sec}. Note that the main difficulty in
proving the second statement is to show that $\lambda _{n}\left( h\right)
\leq \lambda _{2}\left( h\right) $ for every $h\in \left[ 0,1/2\right] $,
which looks evident from looking at Figure \ref{Eigen1234_fig}.

Consequently, in the regime $h\ll 1$, the spectral gap behaves like 
\begin{equation*}
G\left( h\right) \sim \frac{1}{2}h^{2}\text{.}
\end{equation*}

\subsection{Main algebraic structure}

We proceed to describe an\textbf{\ intrinsic} model $\mathbb{W}$ of the
Hermitian vector space $W$, that can be computed as the eigenspace
associated with the maximal eigenvalue of the localized parallel transport
operator $T_{h}$, provided $h\ll 1$. Using this model, the Euclidean inner
products between the viewing directions of every pair of orthonormal frames
can be computed.

\begin{itemize}
\item \textbf{Extrinsic model: }for every point $x\in X$,\textbf{\ }let us
denote by $\delta _{x}:%
%TCIMACRO{\U{2102} }%
%BeginExpansion
\mathbb{C}
%EndExpansion
\rightarrow W$ be the unique complex morphism sending $1\in 
%TCIMACRO{\U{2102} }%
%BeginExpansion
\mathbb{C}
%EndExpansion
$ to the complex vector $f_{1}-if_{2}\in W$.

\item \textbf{Intrinsic model: }we define $\mathbb{W}$ to be the eigenspace
of $T_{h}$ associated with the maximal eigenvalue, which by Theorems \ref%
{spectral_thm} and \ref{gap_thm}, is three-dimensional. For every point $%
x\in X$, there is a map 
\begin{equation*}
\varphi _{x}=\sqrt{2/3}\cdot (ev_{x}|\mathbb{W})^{\ast }:%
%TCIMACRO{\U{2102} }%
%BeginExpansion
\mathbb{C}
%EndExpansion
\rightarrow \mathbb{W}\text{,}
\end{equation*}%
where $ev_{x}:\mathcal{H}\rightarrow 
%TCIMACRO{\U{2102} }%
%BeginExpansion
\mathbb{C}
%EndExpansion
$ is the evaluation morphism at the point $x$, namely, 
\begin{equation*}
ev_{x}\left( f\right) =f\left( x\right) ,
\end{equation*}%
for every $f\in \mathcal{H}$. The pair $\left( \mathbb{W},\left\{ \varphi
_{x}:x\in X\right\} \right) $ is referred to as the intrinsic model of the
vector space $W$.
\end{itemize}

The algebraic structure that underlies the intrinsic classification
algorithm is the canonical morphism 
\begin{equation*}
\tau :W\rightarrow \mathcal{H}\text{,}
\end{equation*}%
defined by%
\begin{equation*}
\tau \left( v\right) \left( x\right) =\sqrt{3/2}\cdot \delta _{x}^{\ast
}\left( v\right) ,
\end{equation*}%
for every $x\in X$. The morphism $\tau $ induces an isomorphism of Hermitian
vector spaces between $W$ equipped with the collection of natural maps $%
\left\{ \delta _{x}:%
%TCIMACRO{\U{2102} }%
%BeginExpansion
\mathbb{C}
%EndExpansion
\rightarrow W\right\} $ and $\mathbb{W}$ equipped with the collection of
maps $\left\{ \varphi _{x}:%
%TCIMACRO{\U{2102} }%
%BeginExpansion
\mathbb{C}
%EndExpansion
\rightarrow \mathbb{W}\right\} $. This is summarized in the following
theorem:

\begin{theorem}
\label{char_thm}The morphism $\tau $ maps $W$ isomorphically, as an
Hermitian vector space, onto the subspace $\mathbb{W}\subset \mathcal{H}$.
Moreover, 
\begin{equation*}
\tau \circ \delta _{x}=\varphi _{x},
\end{equation*}%
for every $x\in X$.
\end{theorem}

For a proof, see Appendix \ref{proofs_sec} (the proof uses the results and
terminology of Section \ref{spectral_sec}).

Using Theorem \ref{char_thm}, we can express in intrinsic terms the inner
product between the viewing directions associated with every ordered pair of
frames. The precise statement is

\begin{theorem}
\label{dist_thm}For every pair of frames $x,y\in X$, we have 
\begin{equation}
\left( \pi \left( x\right) ,\pi \left( y\right) \right) =\left\vert
\left\langle \varphi _{x}\left( v\right) ,\varphi _{y}\left( u\right)
\right\rangle \right\vert -1,  \label{LinAlg2_eq}
\end{equation}%
for any choice of complex numbers $v,u\in 
%TCIMACRO{\U{2102} }%
%BeginExpansion
\mathbb{C}
%EndExpansion
$ of unit norm.
\end{theorem}

For a proof, see Appendix \ref{proofs_sec}. Note that substituting $v=u=1$
in (\ref{LinAlg2_eq}) we obtain (\ref{intrinsic_eq}).

\subsection{Explanation of Theorems \protect\ref{spec_thm} and \protect\ref%
{isom_thm}}

We end this section with an explanation of the two main statements that
appeared in the introduction. The explanation is based on inspecting the
limit when the number of images $N$ goes to infinity. Provided that the
corresponding frames are independently drawn from the normalized Haar
measure on $X$ (empirical assumption); in the limit: the transport matrix $%
\widetilde{T}_{N}$ approaches the localized parallel transport operator $%
T_{h}:\mathcal{H}\rightarrow \mathcal{H}$, for some small value of the
parameter $h$. This implies that the spectral properties of $\widetilde{T}%
_{N}$ for large values of $N$ are governed by the spectral properties of the
operator $T_{h}$ when $h$ lies in the regime $h\ll 1$. In particular,

\begin{itemize}
\item The statement of Theorem \ref{spec_thm} is explained by the fact that
the maximal eigenvalue of $T_{h}$ has multiplicity three (see Theorem \ref%
{spectral_thm}) and that there exists a spectral gap $G(h)\sim h/2$,
separating it from the rest of the spectrum (see Theorem \ref{gap_thm}). The
later property ensures that the numerical computation of this eigenspace
makes sense.

\item The statement of Theorem \ref{isom_thm} is explained by the fact that
the vector space $\mathbb{W}_{N}$ is a numerical approximation of the
theoretical vector space $\mathbb{W}$ and Theorem \ref{char_thm}.
\end{itemize}

\section{Spectral analysis of the localized parallel transport operator\label%
{spectral_sec}}

In this section we study the spectral properties of the localized parallel
transport operator $T_{h}$, mainly focusing on the regime $h\ll 1$. But,
first we need to introduce some preliminaries from representation theory.

\subsection{Isotypic decompositions}

The Hilbert space $\mathcal{H}$, as a unitary representation of the group $%
SO(2)$, admits an isotypic decomposition 
\begin{equation}
\mathcal{H=}\bigoplus\limits_{k\in 
%TCIMACRO{\U{2124} }%
%BeginExpansion
\mathbb{Z}
%EndExpansion
}\mathcal{H}_{k}\text{,}  \label{decomp1_eq}
\end{equation}%
where a function $s\in \mathcal{H}_{k}$ if and only if $s\left(
x\vartriangleleft g\right) =g^{k}s\left( x\right) $, for every $x\in X$ and $%
g\in SO(2)$. In turns, each Hilbert space $\mathcal{H}_{k}$, as a
representation of the group $SO(V)$, admits an isotypic decomposition%
\begin{equation}
\mathcal{H}_{k}=\bigoplus\limits_{n\in 
%TCIMACRO{\U{2115} }%
%BeginExpansion
\mathbb{N}
%EndExpansion
^{\geq 0}}\mathcal{H}_{n,k},  \label{decomp2_eq}
\end{equation}%
where $\mathcal{H}_{n,k}$ denotes the component which is a direct sum of
copies of the unique irreducible representation of $SO(V)$ which is of
dimension $2n+1$. A particularly important property is that each irreducible
representation which appears in (\ref{decomp2_eq}) comes up with
multiplicity one. This is summarized in the following theorem:

\begin{theorem}[Multiplicity one]
\label{mult_prop} If $n<\left \vert k\right \vert $ then $\mathcal{H}%
_{n,k}=0 $. Otherwise, $\mathcal{H}_{n,k}$ is isomorphic to the unique
irreducible representation of $SO(V)$ of dimension $2n+1$.
\end{theorem}

for a proof, see Appendix \ref{proofs_sec}$.$

The following proposition is a direct implication of the equivariance
property of the operator $T_{h}$ and follows from Schur's orthogonality
relations on the group $SO\left( 2\right) $:

\begin{proposition}
\label{kernel_prop}We have 
\begin{equation*}
\bigoplus \limits_{k\neq -1}\mathcal{H}_{k}\subset \ker T_{h}\text{.}
\end{equation*}
\end{proposition}

Consequently, from now on, we will consider $T_{h}$ as an operator from $%
\mathcal{H}_{-1}$ to $\mathcal{H}_{-1}$. Moreover, since for every $n\geq 1$%
, $\mathcal{H}_{n,-1}$ is an irreducible representation of $SO(V)$ and since 
$T_{h}$ commutes with the group action, by Schur's Lemma $T_{h}$ acts on $%
\mathcal{H}_{n,-1}$ as a scalar operator, namely 
\begin{equation*}
T_{h}|\mathcal{H}_{n,-1}=\lambda _{n}\left( h\right) Id\text{.}
\end{equation*}

The reminder of this section is devoted to the computation of\ the
eigenvalues $\lambda _{n}\left( h\right) $. The strategy of the computation
is to choose a point $x_{0}\in X$ and a "good" vector $u_{n}\in \mathcal{H}%
_{n,-1}$ such that $u_{n}\left( x_{0}\right) \neq 0$ and then to use the
relation 
\begin{equation*}
T_{h}\left( u_{n}\right) \left( x_{0}\right) =\lambda _{n}\left( h\right)
u_{n}\left( x_{0}\right) \text{,}
\end{equation*}%
which implies that 
\begin{equation}
\lambda _{n}\left( h\right) =\frac{T_{h}\left( u_{n}\right) \left(
x_{0}\right) }{u_{n}\left( x_{0}\right) }.  \label{eigen_eq}
\end{equation}

\subsection{Set-up}

Fix a frame $x_{0}\in X$, $x_{0}=\left( e_{1},e_{2},e_{3}\right) $. Under
this choice, we can safely identify the group $SO\left( V\right) $ with the
group $SO(3)$ by sending an element $g\in SO(V)$ to the unique element $h\in
SO(3)$ such that $g\vartriangleright x_{0}=x_{0}\vartriangleleft h$. Hence,
from now on, we will consider the frame manifold equipped with a commuting
left and right actions of $SO\left( 3\right) $.

Consider the following elements in the Lie algebra $so\left( 3\right) $: 
\begin{eqnarray*}
A_{1} &=&%
\begin{pmatrix}
0 & 0 & 0 \\ 
0 & 0 & -1 \\ 
0 & 1 & 0%
\end{pmatrix}%
, \\
A_{2} &=&%
\begin{pmatrix}
0 & 0 & 1 \\ 
0 & 0 & 0 \\ 
-1 & 0 & 0%
\end{pmatrix}%
, \\
A_{3} &=&%
\begin{pmatrix}
0 & -1 & 0 \\ 
1 & 0 & 0 \\ 
0 & 0 & 0%
\end{pmatrix}%
.
\end{eqnarray*}

The elements $A_{i}$, $i=1,2,3$ satisfy the relations 
\begin{eqnarray*}
\left[ A_{3},A_{1}\right] &=&A_{2}, \\
\left[ A_{3},A_{2}\right] &=&-A_{1}, \\
\left[ A_{1},A_{2}\right] &=&A_{3}\text{.}
\end{eqnarray*}

Let $\left( H,E,F\right) $ be the following $sl_{2}$ triple in the
complexified Lie algebra $%
%TCIMACRO{\U{2102} }%
%BeginExpansion
\mathbb{C}
%EndExpansion
so\left( 3\right) $:

\begin{eqnarray*}
H &=&-2iA_{3}, \\
E &=&iA_{2}-A_{1}, \\
F &=&A_{1}+iA_{2}\text{.}
\end{eqnarray*}

Finally, let $\left( H^{L},E^{L},F^{L}\right) $ and $\left(
H^{R},E^{R},F^{R}\right) $ be the associated (complexified) vector fields on 
$X$ induced from the left and right action of $SO(3)\ $respectively.

\subsubsection{Spherical coordinates}

We consider the spherical coordinates of the frame manifold $\omega :\left(
0,2\pi \right) \times \left( 0,\pi \right) \times \left( 0,2\pi \right)
\rightarrow X$, given by 
\begin{equation*}
\omega \left( \varphi ,\theta ,\alpha \right) =x_{0}\vartriangleleft
e^{\varphi A_{3}}e^{\theta A_{2}}e^{\alpha A_{3}}\text{.}
\end{equation*}

We have the following formulas

\begin{itemize}
\item The normalized Haar measure on $X$ is given by the density 
\begin{equation*}
\frac{\sin \left( \theta \right) }{2\left( 2\pi \right) ^{2}}d\varphi
d\theta d\alpha \text{.}
\end{equation*}

\item The vector fields $\left( H^{L},E^{L},F^{L}\right) $ are given by 
\begin{eqnarray*}
H^{L} &=&2i\partial _{\varphi }, \\
E^{L} &=&-e^{-i\varphi }\left( i\partial _{\theta }+\cot \left( \theta
\right) \partial _{\varphi }-1/\sin \left( \theta \right) \partial _{\alpha
}\right) , \\
F^{L} &=&-e^{i\varphi }\left( i\partial _{\theta }-\cot \left( \theta
\right) \partial _{\varphi }+1/\sin \left( \theta \right) \partial _{\alpha
}\right) .
\end{eqnarray*}

\item The vector fields $\left( H^{R},E^{R},F^{R}\right) $ are given by 
\begin{eqnarray*}
H^{R} &=&-2i\partial _{\alpha }, \\
E^{R} &=&e^{i\alpha }\left( i\partial _{\theta }+\cot \left( \theta \right)
\partial _{\alpha }-1/\sin \left( \theta \right) \partial _{\varphi }\right)
, \\
F^{R} &=&e^{-i\alpha }\left( i\partial _{\theta }+\cot \left( \theta \right)
\partial _{\alpha }-1/\sin \left( \theta \right) \partial _{\varphi }\right)
.
\end{eqnarray*}
\end{itemize}

\subsection{Choosing a good vector}

\subsubsection{Spherical functions}

Consider the subgroup $T\subset SO\left( 3\right) $ generated by the
infinitesimal element $A_{3}$. For every $k\in 
%TCIMACRO{\U{2124} }%
%BeginExpansion
\mathbb{Z}
%EndExpansion
$ and $n\geq k$, the Hilbert space $\mathcal{H}_{n,k}$ admits an isotypic
decomposition with respect to the left action of $T$: 
\begin{equation*}
\mathcal{H}_{n,k}=\bigoplus\limits_{m=-n}^{n}\mathcal{H}_{n,k}^{m}\text{,}
\end{equation*}%
where a function $s\in \mathcal{H}_{n,k}^{m}$ if and only if $s\left(
e^{-tA_{3}}\vartriangleright x\right) =e^{imt}s\left( x\right) $, for every $%
x\in X$. \ Functions in $\mathcal{H}_{n,k}^{m}$ are usually referred to in
the literature as (generalized) \textit{spherical functions}. Our plan is to
choose for every $n\geq 1$, a spherical function $u_{n}\in \mathcal{H}%
_{n,-1}^{1}$ and exhibit a closed formula for the generating function 
\begin{equation*}
\sum\limits_{n\geq 1}u_{n}t^{n}.
\end{equation*}%
Then, we will use this explicit generating function to compute $u_{n}\left(
x_{0}\right) $ and $T_{h}\left( u_{n}\right) \left( x_{0}\right) $ and use (%
\ref{eigen_eq}) to compute $\lambda _{n}\left( h\right) $.

\subsubsection{Generating function}

For every $n\geq 0$, let $\psi _{n}\in \mathcal{H}_{n,0}^{0}$ be the unique
spherical function such that $\psi _{n}\left( x_{0}\right) =1$. These
functions are the well known spherical harmonics on the sphere. Define the
generating function 
\begin{equation*}
G_{0,0}\left( \varphi ,\theta ,\alpha ,t\right) =\sum\limits_{n\geq 0}\psi
_{n}\left( \varphi ,\theta ,\alpha \right) t^{n}\text{.}
\end{equation*}

The following theorem is taken from \cite{T}:

\begin{theorem}
\label{generating_thm}The function $G_{0,0}$ admits the following formula:%
\begin{equation*}
G_{0,0}\left( \varphi ,\theta ,\alpha ,t\right) =\left( 1-2t\cos \left(
\theta \right) +t^{2}\right) ^{-1/2}\text{.}
\end{equation*}
\end{theorem}

Take $u_{n}=E^{L}F^{R}\psi _{n}$. Note that indeed $u_{n}\in \mathcal{H}%
_{n,-1}^{1}$ and define the generating function 
\begin{equation*}
G_{1,-1}\left( \varphi ,\theta ,\alpha ,t\right) =\sum \limits_{n\geq
1}u_{n}\left( \varphi ,\theta ,\alpha \right) t^{n}\text{.}
\end{equation*}

It follows that, $G_{1,-1}=E^{L}F^{R}G_{0,0}$. Direct calculation, using the
formula in Theorem \ref{generating_thm}, reveals that 
\begin{eqnarray}
G_{1,-1}\left( \varphi ,\theta ,\alpha ,t\right) &=&e^{-i\left( \alpha
+\varphi \right) }[3\sin \left( \theta \right) ^{2}t^{2}\left( 1-2t\cos
\left( \theta \right) +t^{2}\right) ^{-5/2}  \label{generating_eq} \\
&&-t\cos \left( \theta \right) \left( 1-2t\cos \left( \theta \right)
+t^{2}\right) ^{-3/2}  \notag \\
&&-t\left( 1-2t\cos \left( \theta \right) +t^{2}\right) ^{-3/2}]\text{.} 
\notag
\end{eqnarray}

It is enough to consider $G_{1,-1}$ when $\varphi =\alpha =0$. We use the
notation $G_{1,-1}\left( \theta ,t\right) =G_{1,-1}\left( 0,\theta
,0,t\right) $. By (\ref{generating_eq}) 
\begin{eqnarray}
G_{1,-1}\left( \theta ,t\right) &=&3\sin \left( \theta \right)
^{2}t^{2}\left( 1-2t\cos \left( \theta \right) +t^{2}\right) ^{-5/2}
\label{gen1_eq} \\
&&-t\cos \left( \theta \right) \left( 1-2t\cos \left( \theta \right)
+t^{2}\right) ^{-3/2}  \notag \\
&&-t\left( 1-2t\cos \left( \theta \right) +t^{2}\right) ^{-3/2}\text{.} 
\notag
\end{eqnarray}

\subsection{Computation of $u_{n}\left( x_{0}\right) $}

Observe that 
\begin{equation*}
G_{1,-1}\left( 0,t\right) =\sum\limits_{n\geq 1}u_{n}\left( x_{0}\right)
t^{n}\text{.}
\end{equation*}

Direct calculation reveals that 
\begin{eqnarray*}
G_{1,-1}\left( 0,t\right) &=&-2t\left( 1-t\right) ^{-3} \\
&=&-2t\sum \limits_{n\geq 0}\left( \QATOP{-3}{n}\right) \left( -1\right)
^{n}t^{n} \\
&=&-2\sum \limits_{n\geq 1}\left( \QATOP{-3}{n-1}\right) \left( -1\right)
^{n-1}t^{n}.
\end{eqnarray*}

Since $\left( \QATOP{-3}{n-1}\right) =\frac{\left( -1\right) ^{n-1}}{2}%
n\left( n+1\right) $, we obtain 
\begin{equation}
u_{n}\left( x_{0}\right) =-n\left( n+1\right) \text{.}  \label{value_eq}
\end{equation}

\subsection{Computation of $T_{h}\left( u_{n}\right) \left( x_{0}\right) $}

Recall that $h=1-\cos \left( a\right) $.

Using the definition of $T_{h}$, we obtain

\begin{equation*}
T_{h}\left( u_{n}\right) \left( x_{0}\right) =\int\limits_{y\in
B(x_{0},a)}T\left( x_{0},y\right) u_{n}\left( y\right) dy\text{.}
\end{equation*}

Using the spherical coordinates, the integral on the right hand side can be
written as 
\begin{equation*}
\frac{1}{\left( 2\pi \right) ^{2}}\int\limits_{0}^{2\pi }d\varphi
\int\limits_{0}^{a}\frac{\sin \left( \theta \right) }{2}d\theta
\int\limits_{0}^{2\pi }T\left( x_{0},\omega \left( \varphi ,\theta ,\alpha
\right) \right) u_{n}\left( \omega \left( \varphi ,\theta ,\alpha \right)
\right) .
\end{equation*}

First%
\begin{eqnarray}
T\left( x_{0},\omega \left( \varphi ,\theta ,\alpha \right) \right)
&=&T\left( x_{0},x_{0}\vartriangleleft e^{\varphi A_{3}}e^{\theta
A_{2}}e^{\alpha A_{3}}\right)  \label{int1_eq} \\
&=&T\left( x_{0},e^{\varphi A_{3}}\vartriangleright x_{0}\vartriangleleft
e^{\theta A_{2}}e^{\alpha A_{3}}\right)  \notag \\
&=&T(e^{-\varphi A_{3}}\vartriangleright x_{0},x_{0}\vartriangleleft
e^{\theta A_{2}}e^{\alpha A_{3}})  \notag \\
&=&e^{i\varphi }T(x_{0},x_{0}\vartriangleleft e^{\theta A_{2}})e^{i\alpha },
\notag
\end{eqnarray}%
where the third equality uses the invariance property of the transport data
and the second equality uses the equivariance property of the transport data.

Second, since $u_{n}\in \mathcal{H}_{n,-1}^{1}$ we have 
\begin{equation}
u_{n}\left( \omega \left( \varphi ,\theta ,\alpha \right) \right)
=e^{-i\varphi }u_{n}\left( x_{0}\vartriangleleft e^{\theta A_{2}}\right)
e^{-i\alpha }.  \label{int2_eq}
\end{equation}

Combining (\ref{int1_eq}) and (\ref{int2_eq}), we conclude 
\begin{eqnarray}
T_{h}\left( u_{n}\right) \left( x_{0}\right) &=&\int\limits_{0}^{a}\frac{%
\sin \left( \theta \right) }{2}T(x_{0},x_{0}\vartriangleleft e^{\theta
A_{2}})u_{n}\left( x_{0}\vartriangleleft e^{\theta A_{2}}\right) d\theta
\label{int3_eq} \\
&=&\int\limits_{0}^{a}\frac{\sin \left( \theta \right) }{2}u_{n}\left(
x_{0}\vartriangleleft e^{\theta A_{2}}\right) d\theta .  \notag
\end{eqnarray}%
where the second equality uses the fact that $x_{0}\vartriangleleft
e^{\theta A_{2}}$ is the parallel transport of $x_{0}$ along the unique
geodesic connecting $\pi \left( x_{0}\right) $ with $\pi \left(
x_{0}\vartriangleleft e^{\theta A_{2}}\right) $.

Denote 
\begin{equation*}
I_{n}\left( h\right) =\int\limits_{0}^{a}\frac{\sin \left( \theta \right) }{2%
}u_{n}\left( x_{0}\vartriangleleft e^{\theta A_{2}}\right) d\theta \text{.}
\end{equation*}

Define the generating function $I\left( h,t\right) =\sum \limits_{n\geq
0}I_{n}\left( h\right) t^{n}$ and observe that 
\begin{equation*}
I\left( h,t\right) =\int \limits_{0}^{a}\frac{\sin \left( \theta \right) }{2}%
G_{1,-1}\left( \theta ,t\right) d\theta .
\end{equation*}

Direct calculation reveals that 
\begin{eqnarray}
I\left( h,t\right) &=&1/2[h\left( 1+2t\left( h-1\right) +t^{2}\right) ^{-1/2}
\label{integrals_eq} \\
&&-th\left( 2-h\right) \left( 1+2t\left( h-1\right) +t^{2}\right) ^{-3/2} 
\notag \\
&&-t^{-1}(\left( 1+2t\left( h-1\right) +t^{2}\right) ^{1/2}-\left(
1-t\right) )].  \notag
\end{eqnarray}

\subsection{Proof of Theorem \protect\ref{spectral_thm}}

Expanding $I\left( h,t\right) $ with respect to the parameter $t$ reveals
that the function $I_{n}\left( h\right) $ is a polynomial in $h$ of degree $%
n+1$. Then, using Equation (\ref{eigen_eq}), we get%
\begin{equation*}
\lambda _{n}\left( h\right) =-\frac{I_{n}\left( h\right) }{n\left(
n+1\right) }\text{.}
\end{equation*}

In principle, it is possible to obtain a closed formula for $\lambda
_{n}\left( h\right) $ for every $n\geq 1$.

\subsubsection{Quadratic approximation}

We want to compute the first three terms in the Taylor expansion of $\lambda
_{n}\left( h\right) $: 
\begin{equation*}
\lambda _{n}\left( h\right) =\lambda _{n}\left( 0\right) +\partial
_{h}\lambda _{n}\left( 0\right) +\frac{\partial _{h}^{2}\lambda _{n}\left(
0\right) }{2}+O\left( h^{3}\right) \text{.}
\end{equation*}

We have 
\begin{eqnarray*}
\lambda _{n}\left( 0\right) &=&-\frac{I_{n}\left( 0\right) }{n\left(
n+1\right) }, \\
\partial _{h}\lambda _{n}\left( 0\right) &=&-\frac{\partial _{h}I_{n}\left(
0\right) }{n\left( n+1\right) }, \\
\partial _{h}^{2}\lambda _{n}\left( 0\right) &=&-\frac{\partial
_{h}^{2}I_{n}\left( 0\right) }{n\left( n+1\right) }.
\end{eqnarray*}

Observe that 
\begin{equation*}
\partial _{h}^{k}I\left( 0,t\right) =\sum \limits_{n\geq 1}\partial
_{h}^{k}I_{n}\left( 0\right) .
\end{equation*}

Direct computation, using Formula (\ref{integrals_eq}), reveals that 
\begin{eqnarray*}
I\left( 0,t\right) &=&0, \\
\partial _{h}I\left( 0,t\right) &=&-\sum \limits_{n\geq 1}n\left( n+1\right)
t^{n}, \\
\partial _{h}^{2}I\left( 0,t\right) &=&\frac{1}{4}\sum \limits_{n\geq
1}n\left( n+1\right) \left( 1+\left( n+2\right) \left( n-1\right) \right)
t^{n}.
\end{eqnarray*}

Combing all the above yields the desired formula%
\begin{equation*}
\lambda _{n}\left( h\right) =\frac{1}{2}h-\frac{1+\left( n+2\right) \left(
n-1\right) }{8}h^{2}+O\left( h^{3}\right) .
\end{equation*}

This concludes the proof of the theorem.

\appendix

\section{Proofs \label{proofs_sec}}

\subsection{Proof of Theorem \protect\ref{gap_thm}}

The proof is based on two technical lemmas.

\begin{lemma}
\label{gap_pf1_lemma}The following estimates hold:

\begin{enumerate}
\item There exists $h_{1}\in (0,2]$ such that $\lambda _{n}\left( h\right)
\leq \lambda _{1}\left( h\right) $, for every $n\geq 1$ and $h\in \lbrack
0,h_{1}]$.

\item There exists $h_{2}\in (0,2]$ such that $\lambda _{n}\left( h\right)
\leq \lambda _{2}\left( h\right) $, for every $n\geq 2$ and $h\in \left[
0,h_{2}\right] $.
\end{enumerate}
\end{lemma}

The proof appears below.

\begin{lemma}
\label{gap_pf2_lemma}The following estimates hold:

\begin{enumerate}
\item There exists $N_{1}$ such that $\lambda _{n}\left( h\right) \leq
\lambda _{1}\left( h\right) $, for every $n\geq N_{1}$ and $h\in \left[
h_{1},2\right] $.

\item There exists $N_{2}$ such that $\lambda _{n}\left( h\right) \leq
\lambda _{1}\left( h\right) $, for every $n\geq N_{2}$ and $h\in \left[
h_{2},1/2\right] $.
\end{enumerate}
\end{lemma}

The proof appears below.

Granting the validity of these two lemmas we can finish the proof of the
theorem.

First we prove that $\lambda _{n}\left( h\right) \leq \lambda _{1}\left(
h\right) $, for every $n\geq 1$ and $h\in \left[ 0,2\right] $: By Lemmas \ref%
{gap_pf1_lemma},\ref{gap_pf2_lemma} we get that $\lambda _{n}\left( h\right)
\leq \lambda _{1}\left( h\right) $ for every $h\in \left[ 0,2\right] $ when $%
n\geq N_{1}$. Then, we verify directly that $\lambda _{n}\left( h\right)
\leq \lambda _{1}\left( h\right) $ for every $h\in \left[ 0,2\right] $ in
the finitely many cases when $n<N_{1}$.

Similarly, we prove that $\lambda _{n}\left( h\right) \leq \lambda
_{2}\left( h\right) $ for every $n\geq 2$ and $h\in \left[ 0,1/2\right] $:
By Lemmas \ref{gap_pf1_lemma},\ref{gap_pf2_lemma} we get that $\lambda
_{n}\left( h\right) \leq \lambda _{2}\left( h\right) $ for every $h\in \left[
0,1/2\right] $ when $n\geq N_{2}$. Then, we verify directly that $\lambda
_{n}\left( h\right) \leq \lambda _{1}\left( h\right) $ for every $h\in \left[
0,1/2\right] $ in the finitely many cases when $n<N_{2}$.

This concludes the proof of the theorem.

\subsection{Proof of Lemma \protect\ref{gap_pf1_lemma}}

The strategy of the proof is to reduce the statement to known facts about
Legendre polynomials.

Recall $h=1-\cos \left( a\right) $. Here, it will be convenient to consider
the parameter $z=\cos \left( a\right) $, taking values in the interval $%
\left[ -1,1\right] $.

We recall that Legendre polynomials $P_{n}\left( z\right) $, $n\in 
%TCIMACRO{\U{2115} }%
%BeginExpansion
\mathbb{N}
%EndExpansion
$ appear as the coefficients of the generating function 
\begin{equation*}
P(z,t)=\left( t^{2}-2tz+1\right) ^{-(1/2)}.
\end{equation*}

Let 
\begin{equation*}
J_{n}\left( z\right) =\left \{ 
\begin{array}{cc}
\frac{1}{2\left( 1-z\right) } & n=0 \\ 
\frac{1}{2\left( 1-z\right) } & n=1 \\ 
\partial _{z}\lambda _{n-1}\left( z\right) & n\geq 2%
\end{array}%
\right. .
\end{equation*}

Consider the generating function 
\begin{equation*}
J\left( z,t\right) =\sum \limits_{n=0}^{\infty }J_{n}\left( z\right) t^{n}.
\end{equation*}

The function $J\left( z,t\right) $ admits the following closed formula%
\begin{equation}
J\left( z,t\right) =\frac{t+tz+t^{2}+1}{2\left( 1-z\right) }\left(
t^{2}+2tz+1\right) ^{-1/2}\text{.}  \label{gap_pf1_eq}
\end{equation}

Using (\ref{gap_pf1_eq}), we get that for $n\geq 2$%
\begin{equation*}
J_{n}\left( z\right) =\frac{1}{2\left( 1-z\right) }\left( Q_{n}\left(
z\right) +\left( 1+z\right) Q_{n-1}\left( z\right) +Q_{n-2}\left( z\right)
\right) ,
\end{equation*}%
where $Q_{n}\left( z\right) =\left( -1\right) ^{n}P_{n}\left( z\right) $. In
order to prove the lemma, it is enough to show that there exists $z_{0}\in
(-1,1]$ such that for every $z\in \left[ -1,z_{0}\right] $ the following
inequalities hold

\begin{itemize}
\item $Q_{n}\left( z\right) \leq Q_{3}\left( z\right) $ for every $n\geq 3$.

\item $Q_{n}\left( z\right) \leq Q_{2}\left( z\right) $ for every $n\geq 2.$

\item $Q_{n}\left( z\right) \leq Q_{1}\left( z\right) $ for every $n\geq 1.$

\item $Q_{n}\left( z\right) \leq Q_{0}\left( z\right) $ for every $n\geq 0.$
\end{itemize}

These inequalities follow from the following technical proposition.

\begin{proposition}
\label{gap_pf1_prop}Let $n_{0}\in 
%TCIMACRO{\U{2115} }%
%BeginExpansion
\mathbb{N}
%EndExpansion
$. There exists $z_{0}\in (-1,1]$ such that $Q_{n}\left( z\right)
<Q_{n_{0}}\left( z\right) $, for every $z\in \left[ -1,z_{0}\right] $ and $%
n\geq n_{0}$.
\end{proposition}

The proof appears below.

Take $h_{0}=h_{1}=1+z_{0}$. Granting Proposition \ref{gap_pf1_eq}, verify
that $J_{n}\left( z\right) \leq J_{2}\left( z\right) $, for $n\geq 2$, $z\in %
\left[ -1,z_{0}\right] $ which implies that $\lambda _{n}\left( h\right)
\leq \lambda _{1}\left( h\right) $, for $n\geq 2$, $h\in \left[ 0,h_{0}%
\right] $ and $J_{n}\left( z\right) \leq J_{3}\left( z\right) $, for $n\geq
3 $, $z\in \left[ -1,z_{0}\right] $ which implies that $\lambda _{n}\left(
h\right) \leq \lambda _{2}\left( h\right) $, for $n\geq 3$, $h\in \left[
0,h_{0}\right] $.

This concludes the proof of the Lemma.

\subsubsection{Proof of Proposition \protect\ref{gap_pf1_prop}}

Denote by $a_{1}<a_{2}<...<a_{n}$ the zeroes of $Q_{n}\left( \cos \left(
a\right) \right) $ and by $\mu _{1}<\mu _{2}<...<\mu _{n-1}$ the local
extrema of $Q_{n}\left( \cos \left( a\right) \right) $.

The following properties of the polynomials\ $Q_{n}$ are implied from known
facts about Legendre polynomials (Properties 1 and 2 can be verified
directly), that can be found for example in the book \cite{SZ}:

\begin{description}
\item[Property 1] $a_{i}<\mu _{i}<a_{i+1}$, for $i=1,..,n-1$.

\item[Property 2] $Q_{n}\left( -1\right) =1$ and $\partial _{z}Q_{n+1}\left(
-1\right) <\partial _{z}Q_{n}\left( -1\right) <0$, for $n\in 
%TCIMACRO{\U{2115} }%
%BeginExpansion
\mathbb{N}
%EndExpansion
$.

\item[Property 3] $\left \vert Q_{n}\left( \cos \left( \mu _{i}\right)
\right) \right \vert \geq \left \vert Q_{n}\left( \cos \left( \mu
_{i+1}\right) \right) \right \vert $, for $i=1,..,[n/2]$.

\item[Property 4] $\left( i-1/2\right) \pi /n\leq a_{i}\leq i\pi /(n+1)$,
for $i=1,..,\left[ n/2\right] $.

\item[Property 5] $\sin \left( a\right) ^{1/2}\cdot \left \vert Q_{n}\left(
\cos \left( a\right) \right) \right \vert <\sqrt{2/\pi n}$, for $a\in \left[
0,\pi \right] $.
\end{description}

Granting these facts, we can finish the proof.

By Properties 1,4%
\begin{equation*}
\frac{\pi }{2n}<\mu _{1}<\frac{2\pi }{n+1}\text{.}
\end{equation*}

We assume that $n$ is large enough so that, for some small $\epsilon >0$ 
\begin{equation*}
\sin \left( a_{1}\right) \geq \left( 1-\epsilon \right) a_{1},
\end{equation*}

In particular, this is the situation when $n_{0}\geq N$, for some fixed $%
N=N_{\epsilon }$.

By Property 5 
\begin{eqnarray*}
\left \vert Q_{n}\left( \cos \left( \mu _{1}\right) \right) \right \vert &<&%
\sqrt{2/\pi n}\cdot \sin \left( \mu _{1}\right) ^{-1/2} \\
&<&\sqrt{2/\pi n}\cdot \sin \left( a_{1}\right) ^{-1/2} \\
&<&\sqrt{2/\pi n}\cdot \left( \left( 1-\epsilon \right) a_{1}\right) ^{-1/2}=%
\frac{2}{\pi \sqrt{1-\epsilon }}\text{.}
\end{eqnarray*}

Let $a_{0}\in \left( 0,\pi \right) $ be such that $Q_{n_{0}}\left( \cos
\left( a\right) \right) >2/\pi \sqrt{1-\epsilon }$, for every $a<a_{0}$.
Take $z_{0}=\cos \left( a_{0}\right) $.

Finally, in the finitely many cases where $n_{0}\leq n\leq N$, the
inequality $Q_{n}\left( z\right) <Q_{n_{0}}\left( z\right) $ can be verified
directly.

This concludes the proof of the proposition.

\subsection{Proof of Lemma \protect\ref{gap_pf2_lemma}}

We have the following identity: 
\begin{equation}
tr\left( T_{h}^{2}\right) =\frac{h}{2},  \label{trace_eq}
\end{equation}

for every $h\in \left[ 0,2\right] $. The proof of (\ref{trace_eq}) is by
direct calculation:%
\begin{eqnarray*}
tr\left( T_{h}^{2}\right) &=&\int\limits_{x\in X}T_{h}^{2}\left( x,x\right)
dx= \\
&=&\int\limits_{x\in X}\mu _{Haar}\int\limits_{y\in B\left( x,a\right)
}T_{h}\left( x,y\right) \circ T_{h}\left( y,x\right) dx.
\end{eqnarray*}

Since $T_{h}\left( x,y\right) =T_{h}\left( y,x\right) ^{-1}$ (symmetry
property), we get 
\begin{equation*}
tr\left( T_{h}^{2}\right) =\int\limits_{x\in X}\int\limits_{y\in B\left(
x,a\right) }dxdy=\int\limits_{0}^{a}\frac{\sin \left( \theta \right) }{2}%
d\theta =\frac{1-\cos \left( a\right) }{2}.
\end{equation*}

Substituting, $a=\cos ^{-1}\left( 1-h\right) $, we get the desired formula $%
tr\left( T_{h}^{2}\right) =h/2$.

On the other hand, 
\begin{equation}
tr\left( T_{h}^{2}\right) =\sum \limits_{n=1}^{\infty }tr\left( T_{h|%
\mathcal{H}_{n,-1}}^{2}\right) =\sum \limits_{n=1}^{\infty }\left(
2n+1\right) \lambda _{n}\left( h\right) ^{2}\text{.}  \label{trace1_eq}
\end{equation}

From (\ref{trace_eq}) and (\ref{trace1_eq}) we obtain the following upper
bound 
\begin{equation}
\lambda _{n}\left( h\right) \leq \frac{\sqrt{h}}{\sqrt{4n+2}}\text{.}
\label{upper_eq}
\end{equation}

Now we can finish the proof.

First estimate: We know that $\lambda _{1}\left( h\right) =h/2-h^{2}/8$,
hence, one can verify directly that there exists $N_{1}$ such that $\sqrt{h}/%
\sqrt{4n+2}\leq \lambda _{1}\left( h\right) $ for every $n\geq N_{1}$ and $%
h\in \left[ h_{1},2\right] $, which implies by (\ref{upper_eq}) that $%
\lambda _{n}\left( h\right) \leq \lambda _{1}\left( h\right) $ for every $%
n\geq N_{1}$ and $h\in \left[ h_{1},2\right] $.

Second estimate: We know that $\lambda _{2}\left( h\right)
=h/2-5h^{2}/8+h^{3}/6$, therefore, one can verify directly that there exists 
$N_{2}$ such that $\sqrt{h}/\sqrt{4n+2}\leq \lambda _{2}\left( h\right) $
for every $n\geq N_{2}$ and $h\in \left[ h_{2},1/2\right] $, which implies
by (\ref{upper_eq}) that $\lambda _{n}\left( h\right) \leq \lambda
_{2}\left( h\right) $ for every $n\geq N_{2}$ and $h\in \left[ h_{2},1/2%
\right] $.

This concludes the proof of the Lemma.

\subsection{Proof of Theorem \protect\ref{char_thm}}

We begin by proving that $\tau $ maps $W=%
%TCIMACRO{\U{2102} }%
%BeginExpansion
\mathbb{C}
%EndExpansion
V$ isomorphically, as an Hermitian space, onto $\mathbb{W}=\mathcal{H}\left(
\lambda _{\max }\left( h\right) \right) $.

The crucial observation is, that $\mathcal{H}\left( \lambda _{\max }\left(
h\right) \right) $ coincide with the isotypic subspace $\mathcal{H}_{1,-1}$
(see Section \ref{spectral_sec}). Consider the morphism $\alpha =\sqrt{2/3}%
\cdot \tau :W\rightarrow \mathcal{H}$, given by 
\begin{equation*}
\alpha \left( v\right) \left( x\right) =\delta _{x}^{\ast }\left( v\right) 
\text{.}
\end{equation*}

First claim is, that $\func{Im}\alpha \subset \mathcal{H}_{-1}$, namely,
that $\delta _{x\vartriangleleft g}^{\ast }\left( v\right) =g^{-1}\delta
_{x}^{\ast }\left( v\right) $, for every $v\in W$, $x\in X$ and $g\in SO(2)$%
. Denote by $\left\langle \cdot ,\cdot \right\rangle _{std}$ the standard
Hermitian product on $%
%TCIMACRO{\U{2102} }%
%BeginExpansion
\mathbb{C}
%EndExpansion
$. Now write 
\begin{eqnarray*}
\langle \delta _{x\vartriangleleft g}^{\ast }\left( v\right) ,z\rangle
_{std} &=&\left\langle v,\delta _{x\vartriangleleft g}\left( z\right)
\right\rangle =\left\langle v,\delta _{x}\left( gz\right) \right\rangle \\
&=&\langle \delta _{x}^{\ast }\left( v\right) ,gz\rangle _{std}=\langle
g^{-1}\delta _{x}^{\ast }\left( v\right) ,z\rangle _{std}\text{.}
\end{eqnarray*}

Second claim is, that $\alpha $ is a morphism of $SO\left( V\right) $
representations, namely, that $\delta _{x}^{\ast }\left( gv\right) =\delta
_{g^{-1}\vartriangleright x}\left( v\right) $, for every $v\in W$, $x\in X$
and $g\in SO(V)$. This statement follows from 
\begin{eqnarray*}
\left\langle \delta _{x}^{\ast }\left( gv\right) ,z\right\rangle _{std}
&=&\left\langle gv,\delta _{x}\left( z\right) \right\rangle =\left\langle
v,g^{-1}\delta _{x}\left( z\right) \right\rangle \\
&=&\langle v,\delta _{g^{-1}\vartriangleright x}\left( z\right) \rangle
=\langle \delta _{g^{-1}\vartriangleright x}^{\ast }\left( v\right)
,z\rangle _{std}.
\end{eqnarray*}

Consequently, the morphism $\alpha $ maps $W$ isomorphically, as a unitary
representation of $SO\left( V\right) $, onto $\mathcal{H}_{1,-1}$, which is
the unique copy of the three-dimensional representation of $SO\left(
V\right) $ in $\mathcal{H}_{-1}$. In turns, this implies that, up to a
scalar, $\alpha $ and, hence $\tau $, are isomorphisms of Hermitian spaces.
In order to complete the proof it is enough to show that 
\begin{equation*}
tr\left( \tau ^{\ast }\circ \tau \right) =3\text{.}
\end{equation*}

This follows from 
\begin{eqnarray*}
tr\left( \tau \circ \tau ^{\ast }\right) &=&\frac{3}{2}tr\left( \alpha
^{\ast }\circ \alpha \right) \\
&=&\frac{3}{2}\int\limits_{v\in S\left( W\right) }\left\langle \alpha ^{\ast
}\circ \alpha \left( v\right) ,v\right\rangle _{\mathcal{H}}dv \\
&=&\frac{3}{2}\int\limits_{v\in S\left( W\right) }\left\langle \alpha \left(
v\right) ,\alpha \left( v\right) \right\rangle _{\mathcal{H}}dv \\
&=&\frac{3}{2}\int\limits_{v\in S\left( W\right) }\int\limits_{x\in
X}\left\langle \delta _{x}^{\ast }\left( v\right) ,\delta _{x}^{\ast }\left(
v\right) \right\rangle _{std}dvdx \\
&=&\frac{3}{2}\int\limits_{v\in S\left( W\right) }\int\limits_{x\in
X}2dvdx=3.
\end{eqnarray*}%
where $dv$ denotes the normalized Haar measure on the five dimensional
sphere $S\left( W\right) $.

Next, we prove that $\tau \circ \delta _{x}=\varphi _{x}$, for every $x\in X$%
. The starting point is the equation $ev_{x}|W\circ \alpha =\delta
_{x}^{\ast }$, which follows from the definition of the morphism $\alpha $
and the fact that $\func{Im}\alpha =\mathbb{W}$. This implies that $\varphi
_{x}^{\ast }\circ \tau =\delta _{x}^{\ast }$. The statement now follows from%
\begin{eqnarray*}
\varphi _{x}^{\ast }\circ \tau &=&\delta _{x}^{\ast }\Rightarrow \varphi
_{x}^{\ast }\circ \left( \tau \circ \tau ^{\ast }\right) =\delta _{x}^{\ast
}\circ \tau ^{\ast } \\
&\Rightarrow &\varphi _{x}^{\ast }=\delta _{x}^{\ast }\circ \tau ^{\ast
}\Rightarrow \varphi _{x}=\tau \circ \delta _{x}\text{.}
\end{eqnarray*}

This concludes the proof of the theorem.

\subsection{Proof of Theorem \protect\ref{dist_thm}}

We use the following terminology: for every $x\in X$, $x=\left(
e_{1},e_{2},e_{3}\right) $, we denote by $\widetilde{\delta }_{x}:%
%TCIMACRO{\U{2102} }%
%BeginExpansion
\mathbb{C}
%EndExpansion
\rightarrow V$ the map given by $\widetilde{\delta }_{x}\left( p+iq\right)
=pe_{1}+qe_{2}$. We observe that $\delta _{x}\left( v\right) =$ $\widetilde{%
\delta }_{x}\left( v\right) -i\widetilde{\delta }_{x}\left( iv\right) $, for
every $v\in 
%TCIMACRO{\U{2102} }%
%BeginExpansion
\mathbb{C}
%EndExpansion
$.

We proceed with the proof. Let $x,y\in X$. Choose unit vectors $%
v_{x},v_{y}\in 
%TCIMACRO{\U{2102} }%
%BeginExpansion
\mathbb{C}
%EndExpansion
$ such that $\widetilde{\delta }_{x}\left( v_{x}\right) =$ $\widetilde{%
\delta }_{y}\left( v_{y}\right) =v$.

Write%
\begin{eqnarray}
\left\langle \delta _{x}\left( v_{x}\right) ,\delta _{y}\left( v_{y}\right)
\right\rangle &=&\left\langle \widetilde{\delta }_{x}\left( v_{x}\right) -i%
\widetilde{\delta }_{x}\left( iv_{x}\right) ,\widetilde{\delta }_{y}\left(
v_{y}\right) -i\widetilde{\delta }_{y}\left( iv_{y}\right) \right\rangle 
\notag \\
&=&\left( \widetilde{\delta }_{x}\left( v_{x}\right) ,\widetilde{\delta }%
_{y}\left( v_{y}\right) \right) +\left( \widetilde{\delta }_{x}\left(
iv_{x}\right) ,\widetilde{\delta }_{y}\left( iv_{y}\right) \right)
\label{dist1_eq} \\
&&-i\left( \widetilde{\delta }_{x}\left( iv_{x}\right) ,\widetilde{\delta }%
_{y}\left( v_{y}\right) \right) +i\left( \widetilde{\delta }_{x}\left(
v_{x}\right) ,\widetilde{\delta }_{y}\left( iv_{y}\right) \right) .  \notag
\end{eqnarray}

For every frame $z\in X$ and vector\ $v_{z}\in 
%TCIMACRO{\U{2102} }%
%BeginExpansion
\mathbb{C}
%EndExpansion
$, the following identity can be easily verified: 
\begin{equation*}
\widetilde{\delta }_{z}\left( iv_{z}\right) =\pi \left( z\right) \times 
\widetilde{\delta }_{z}\left( v_{z}\right) .
\end{equation*}

This implies that 
\begin{eqnarray*}
\widetilde{\delta }_{x}\left( iv_{x}\right) &=&\pi \left( x\right) \times 
\widetilde{\delta }_{x}\left( v_{x}\right) =\pi \left( x\right) \times v, \\
\widetilde{\delta }_{y}\left( iv_{y}\right) &=&\pi \left( y\right) \times 
\widetilde{\delta }_{y}\left( v_{y}\right) =\pi \left( y\right) \times v.
\end{eqnarray*}

Combining these identities with Equation (\ref{dist1_eq}), we obtain%
\begin{equation*}
\left\langle \delta _{x}\left( v_{x}\right) ,\delta _{y}\left( v_{y}\right)
\right\rangle =\left( v,v\right) +\left( \pi \left( x\right) \times v,\pi
\left( y\right) \times v\right) -i\left( \pi \left( x\right) \times
v,v\right) +i\left( v,\pi \left( y\right) \times v\right) .
\end{equation*}

Since $v\in \func{Im}\widetilde{\delta }_{x}\cap \func{Im}\widetilde{\delta }%
_{y}$, it follows that $\left( \pi \left( x\right) \times v,v\right) =\left(
v,\pi \left( y\right) \times v\right) =0$. In addition, 
\begin{equation*}
\left( \pi \left( x\right) \times v,\pi \left( y\right) \times v\right)
=\det 
\begin{pmatrix}
\left( \pi \left( x\right) ,\pi \left( y\right) \right) & \left( \pi \left(
x\right) ,v\right) \\ 
\left( \pi \left( y\right) ,v\right) & \left( v,v\right)%
\end{pmatrix}%
=\left( \pi \left( x\right) ,\pi \left( y\right) \right) \text{.}
\end{equation*}

Thus, we obtain that $\left\langle \delta _{x}\left( v_{x}\right) ,\delta
_{y}\left( v_{y}\right) \right\rangle =1+\left( \pi \left( x\right) ,\pi
\left( y\right) \right) $. Since the right hand side is always $\geq 0$ it
follows that 
\begin{equation}
\left\vert \left\langle \delta _{x}\left( v_{x}\right) ,\delta _{y}\left(
v_{y}\right) \right\rangle \right\vert =1+\left( \pi \left( x\right) ,\pi
\left( y\right) \right) .  \label{dist2_eq}
\end{equation}

Now, notice that the left hand side of \ref{dist2_eq} does not depend on the
choice of the unit vectors $v_{x}$ and $v_{y}$.

To finish the proof, we use the isomorphism $\tau $ which satisfies $\tau
\circ \delta _{x}=\varphi _{x}$ for every $x\in X$, and get%
\begin{equation*}
\left\vert \left\langle \varphi _{x}\left( v_{x}\right) ,\varphi _{y}\left(
v_{y}\right) \right\rangle \right\vert =1+\left( \pi \left( x\right) ,\pi
\left( y\right) \right) .
\end{equation*}

This concludes the proof of the theorem.

\subsection{Proof of Proposition \protect\ref{mult_prop}}

The basic observation is, that $\mathcal{H}$, as a representation of $%
SO(V)\times SO(3)$, admits the following isotypic decomposition%
\begin{equation*}
\mathcal{H=}\bigoplus \limits_{n=0}^{\infty }V_{n}\otimes U_{n}\text{,}
\end{equation*}%
where $V_{n}$ is the unique irreducible representation of $SO(V)$ of
dimension $2n+1$, and, similarly, $U_{n}$ is the unique irreducible
representation of $SO(3)$ of dimension $2n+1$. This assertion, principally,
follows from the Peter Weyl Theorem for the regular representation of $SO(3)$%
.

This implies that the isotypic decomposition of $\mathcal{H}_{k}$ takes the
following form%
\begin{equation*}
\mathcal{H}_{k}=\bigoplus \limits_{n=0}^{\infty }V_{n}\otimes U_{n}^{k}\text{%
,}
\end{equation*}%
where $U_{n}^{k}$ is the weight $k$ space with respect to the action $%
SO(2)\subset SO(3)$. The statement now follows from\ the following standard
fact about the weight decomposition: 
\begin{equation*}
\dim U_{n}^{k}=\left \{ 
\begin{array}{cc}
0 & n<k \\ 
1 & n\geq k%
\end{array}%
\right. .
\end{equation*}

This concludes the proof of the theorem.

\bigskip

\end{document}